\documentclass{amsart}
\usepackage{delimset, amssymb, enumitem, caption, mathtools, thmtools}
\setlist{itemsep=4pt, topsep=0pt, leftmargin=17pt, listparindent=11pt}
\newcommand{\fo}{\operatorname{fo}}
\newcommand{\lo}{\operatorname{lo}}

\usepackage{xcolor}
\definecolor{PKU}{cmyk}{0, 1, 1, .45}
\definecolor{BIT}{cmyk}{1, 0, 1, 0}
\usepackage[colorlinks, citecolor=BIT, linkcolor=PKU]{hyperref}
\usepackage{float}
    
\usepackage{tikz}
\usetikzlibrary{decorations.pathreplacing, calc, positioning, arrows.meta}
\tikzset{
edge/.style={semithick},
ball/.style={shape=circle, minimum size=1mm, ball color=black, inner sep=0.55},
ellipsis/.style={shape=circle, fill, inner sep=.5}}
\def\r{1}
\def\eps{.2}

\usepackage[capitalize]{cleveref}
\crefname{ineq}{Ineq.}{Ineqs.}
\Crefname{ineq}{Inequality}{Inequalities}
\creflabelformat{ineq}{~\upshape(#2#1#3)}
\crefname{def}{Def.}{Defs.}
\Crefname{def}{Definition}{Definitions}
\creflabelformat{def}{~\upshape(#2#1#3)}

\makeatletter
\newcommand\creflabel[2][\@currentcounter]{%
 \crefalias{\@currentcounter}{#1}\label{#2}}
\makeatother 

\crefname{itm}{}{}
\Crefname{itm}{Item}{Items}
\creflabelformat{itm}{~\upshape(#2#1#3)}

\crefname{conjecture}{Conjecture}{Conjectures}
\newtheorem{theorem}{Theorem}[section]
\newtheorem{lemma}[theorem]{Lemma}
\newtheorem{conjecture}[theorem]{Conjecture}

\numberwithin{equation}{section}
\numberwithin{figure}{section}

\allowbreak
\allowdisplaybreaks

\usepackage[numbers, sort&compress, nonamebreak, merge, elide, longnamesfirst]{natbib}

\author[D.Q.B.~Tang]{Davion Q.B. Tang}
\address[Davion Q.B. Tang]{School of Mathematics and Statistics, Beijing Institute of Technology, Beijing 102400, P. R. China}
\email{davion@bit.edu.cn}

\author[D.G.L.~Wang]{David G.L. Wang}
\address[David G.L. Wang]{School of Mathematics and Statistics \& Beijing Key Laboratory on MCAACI, Beijing Institute of Technology; MIIT Key Laboratory of Mathematical Theory and Computation in Information Security, Beijing 102400, P. R. China}
\email{glw@bit.edu.cn}

\author[M.M.Y. Wang]{Monica M.Y. Wang}
\address[Monica M.Y. Wang]{School of Basic Education, Beijing Institute of Graphic Communication, Beijing 102600, P. R. China}
\email{myw@bigc.edu.cn}

\keywords{chromatic symmetric function; $e$-positivity; Stanley and Stembridge's conjecture; the composition method}

\subjclass[2020]{05E05, 05C15, 05A20}

\thanks{David Wang is the corresponding author, and is supported by the General Program of National Natural Science Foundation of China (Grant No.~12171034). Monica Wang is supported by the Funding Project of Beijing Institute of Graphic Communication (No.~Ea202512).}

\title[]{The spiders $S(4m+2,\,2m,\,1)$ are $e$-positive}

\begin{document}
\bibliographystyle{abbrvnat}
\begin{abstract}
Using the composition method, we establish the $e$-positivity of spiders of the form $S(4m+2,\, 2m,\, 1)$, as conjectured by Aliniaeifard, van Willigenburg and Wang. Following a divide-and-conquer strategy, we group one or two positive-coefficient $e_J$-terms with each negative-coefficient $e_I$-term, where each $J$ is a composition obtained by rearranging the parts of the composition $I$, and prove that the resulting total coefficient is nonnegative. Our main contribution is the explicit construction of the required injections.
\end{abstract}
\maketitle


\section{Introduction}

In his seminal paper, \citet{Sta95} introduced the chromatic symmetric function $X_G$ of a graph $G$ as a symmetric-function analogue of Birkhoff's chromatic polynomial. From the beginning, a basic problem has been to understand the expansion of $X_G$ in the standard bases of the algebra of symmetric functions. Given a basis $\{b_{\lambda}\}$ for symmetric functions, a symmetric function $f$ is said to be $b$-positive if every $b$-coefficient of $f$ is nonnegative. A graph is said to be $b$-positive if its chromatic symmetric function is $b$-positive. \citet[Section 5]{Sta95} asked the question: Which graphs are $e$-positive?

The $e$-basis is especially central because of the conjecture of \citet{SS93}. In the language of chromatic symmetric functions, Stanley and Stembridge's conjecture asserts that the incomparability graph of any $(3+1)$-free poset is $e$-positive. Although this is a rather special class of graphs, the conjecture remained open for more than three decades, and was recently confirmed by \citet[Corollary~1.8]{Hik24X}. This illustrates the intrinsic difficulty of proving $e$-positivity.

A complete classification of $e$-positive graphs appears to be out of reach at present. Consequently, much of the current work focuses on identifying families of graphs whose chromatic symmetric functions are $e$-positive. For example, \citet{TV25X} established the $e$-positivity of the sum graph of natural unit interval graphs and cycles, and that of noncrossing cycle-chord graphs. Very recently, \citet{CHW26} showed the $e$-positivity of clocks, which do not belong to any of the graph classes mentioned above. The general question of which graphs are $e$-positive is still open.

Trees form a particularly important testing ground for $e$-positivity. On the one hand, they are among the structurally simplest connected graphs. On the other hand, $e$-positivity appears to be rare among trees, and even finding an infinite family of $e$-positive trees is difficult. This scarcity is already visible from direct computation: among the $11{,}989{,}764$ non-isomorphic connected graphs of order at most~$10$, $3{,}396{,}886$ are $e$-positive; whereas among the $522{,}959$ non-isomorphic trees of order at most~$19$, only $51$ are $e$-positive. In these ranges, the corresponding proportions of $e$-positive graphs are
\[
\frac{3{,}396{,}886}{11{,}989{,}764}\approx 28.3316\%
\quad\text{and}\quad
\frac{51}{522{,}959}\approx 0.009752\%,
\]
respectively. These comparisons help explain the interest in explicit infinite families of $e$-positive trees, including the family considered in this paper.

A parallel direction is the study of non-$e$-positivity of trees. \citet[Conjecture~6.1]{DSv20} conjectured that every $e$-positive tree $T$ has maximum degree $\Delta(T)$ at most~$3$. \citet[Theorems~1.5 and~1.6]{Tom24X} made significant progress by proving non-$e$-positivity for the following cases:
\begin{itemize}
\item trees with $\Delta(T)\ge 5$;
\item trees with $\Delta(T)=4$ and at least one degree-$4$ vertex that is not adjacent to any leaf;
\item spiders with exactly four legs.
\end{itemize}
Thus, except for the settled case of $4$-legged spiders, the remaining maximum-degree-$4$ trees are those in which every degree-$4$ vertex is adjacent to a leaf.

Indeed, spiders play an essential role in this direction. For any partition $\lambda=\lambda_1\dotsm\lambda_d$ of $n-1$ with $d\ge 3$, the \emph{spider}~$S(\lambda)$ is the $n$-vertex tree consisting of paths of lengths $\lambda_1,\dots,\lambda_d$ with a common end. For instance, the spider $S(4m+2,\,2m,\,1)$ is obtained by adding a pendant vertex to the path of order $6m+3$, see \cref{fig:S}. The importance of spiders is explained by connected partition obstructions. \citet[Lemma~13]{DSv20} showed that, after deleting a vertex of degree at least~$3$ from a connected graph $G$ and recording the sizes of the resulting connected components as a partition~$\lambda$, connected partition obstructions for the spider $S(\lambda)$ give corresponding obstructions for the $e$-positivity of~$G$. Together with the criterion of \citet[Proposition~1.3.3]{Wol97D}, which states that any connected $e$-positive graph has a connected partition of every type, this makes spiders a basic testing class for $e$-positivity questions on trees.

\begin{figure}[htbp]
\begin{tikzpicture}
\node (O) at (0,0) [ball]{};

\node (a1) at (-\r, 0) [ball]{};
\node (b1) at (\r, 0) [ball]{};
\node (c1) at (0, \r) [ball]{};
\node (a2) at (-\r*2, 0) [ball]{};
\node (a4) at (-\r*4, 0) [ball]{};

\node (a3) at (-\r*3, 0) [ellipsis]{};
\node[below=5pt] at (a3) {length $4m+2$};

\node (b2) at (\r*2, 0) [ellipsis]{};
\node[below=5pt] at (b2) {length $2m$};

\node (a3l) at (-\r*3+\eps, 0) [ellipsis]{};
\node (b2l) at (\r*2-\eps, 0) [ellipsis]{};

\node (a3r) at (-\r*3-\eps, 0) [ellipsis]{};
\node (b2r) at (\r*2+\eps, 0) [ellipsis]{};

\node (a5) at (-\r*5, 0) [ball]{};
\node (b3) at (\r*3, 0) [ball]{};

\draw[edge] (O) -- (a1) -- (a2) -- ($(a3l) - (-\eps, 0)$);
\draw[edge] (O) -- (b1) -- ($(b2l) + (-\eps, 0)$);
\draw[edge] (O) -- (c1);

\draw[edge] ($(a3r) - (\eps, 0)$) -- (a4) --(a5);
\draw[edge] ($(b2r) + (\eps, 0)$) -- (b3);
\end{tikzpicture}
\caption{The spider $S(4m+2,\,2m,\,1)$.}
\label{fig:S}
\end{figure}

\citet[Conjecture~6.1]{WW23-DAM} conjectured that a $3$-legged spider $S(a,\,b,\,c)$ with $a\ge b\ge c\ge 3$ is $e$-positive if and only if it is $S(8,5,3)$ or $S(14,9,5)$. In contrast, for spiders $S(a,\,b,\,c)$ with $c\in\{1,2\}$, the existing results suggest a much less uniform picture; there appears to be no clear pattern governing whether a given such spider is $e$-positive. Against this background, \citet[Conjecture~6.2]{Zhe22} presented \cref{conj:epos:4m+2.2m.1} and attributed it to Aliniaeifard, van Willigenburg and Wang.

\begin{conjecture}[Aliniaeifard, van Willigenburg and Wang]\label{conj:epos:4m+2.2m.1}
For any integer $m\ge 1$, the spider $S(4m+2,\,2m,\,1)$ is $e$-positive.
\end{conjecture}

For example, the $e$-positivity of $S(6,\,2,\,1)$ is clear from the $e$-expansion
\begin{align*}
X_{S(6,2,1)}
&=10e_{10}
+17e_{91}
+22e_{82}
+7e_{81^2}
+11e_{73}
+24e_{721}
+38e_{64}
+32e_{631}
+26e_{62^2}
+5e_{621^2}\\
&\quad
+20e_{5^2}
+55e_{541}
+37e_{532}
+16e_{531^2}
+20e_{52^21}
+42e_{4^22}
+9e_{4^21^2}
+e_{43^2}
+59e_{4321}\\
&\quad
+22e_{42^3}
+3e_{42^21^2}
+8e_{3^31}
+9e_{3^22^2}
+8e_{3^221^2}
+9e_{32^31}
+2e_{2^5}.
\end{align*}

The following more general conjecture can be found in \citet[Conjecture~6.3]{Zhe22}.

\begin{conjecture}[Aliniaeifard, van Willigenburg and Wang]\label{conj2}
For any integers $m,n\ge 1$, the spider $S(n(n!m+1),\,n!m,\,1)$ is $e$-positive.
\end{conjecture}
It is true for $n=1$, since
\[
X_{S(m+1,\,m,\,1)}
=\frac{1}{2}
\brk1{
X_{B_{(m+1)(m+1)}}
+e_1 
X_{P_{2m+2}}}.
\]
Here $B_{(m+1)(m+1)}$ denotes the generalized bull graph obtained by adding
a path of length $m$ to each of two distinct vertices of a triangle. This graph
is $e$-positive; see \citet[Page~2681]{DFv20}. For $n=2$, \cref{conj2} reduces to \cref{conj:epos:4m+2.2m.1}.

This paper confirms \cref{conj:epos:4m+2.2m.1}. Our proof uses the composition method and proceeds by decomposing the chromatic symmetric function of $S(4m+2,\,2m,\,1)$ into several components. The first components are directly $e$-positive, and the remaining component is handled by a divide-and-conquer argument on its $e_I$-terms. More precisely, we group one or two positive $e_J$-terms with each negative $e_I$-term, where the composition $J$ is obtained by rearranging the parts of $I$, and then prove that the total coefficient is nonnegative. The main point is an explicit construction of the required injection, which shows that the positive terms used in this cancellation procedure are all distinct.

In \cref{sec:preliminary}, we recall the necessary notation and the formulas for chromatic symmetric functions of paths and $3$-legged spiders. The proof of \cref{conj:epos:4m+2.2m.1} is given in \cref{sec:proof}.

\section{Preliminary}\label{sec:preliminary}
Let $n$ be a nonnegative integer. A \emph{composition} of $n$ is a sequence of positive integers that sum to~$n$, denoted 
\[
I=i_1\cdots i_l\vDash n.
\]
In particular, the unique composition of $0$ is the empty sequence. We denote 
\[
\abs{I}=n
\quad\text{and}\quad
\ell(I)=l.
\]
The integers $i_k$ are called \emph{parts} of~$I$, and we denote the $j$th last part $i_{\ell(I)+1-j}$ by $i_{-j}$ for notational convenience. Whenever a capital letter such as $I,J\dots$ denotes a composition, we use the corresponding small letter counterparts $i,j\dots$ with integer subscripts to denote its parts. Following \citet[Eq.~(2.12)]{WZ25}, the \emph{$a$-surplus} of $I$ is the number
\[
\Theta_I(a)
=
\min\brk[c]1{
\abs{i_1\dotsm i_k}\colon
0\le k\le \ell(I),\
\abs{i_1\dotsm i_k}\ge a}
-a.
\]
It is clear that $\Theta_I(a)\ge 0$. A \emph{partition} of $n$ is a composition of $n$ in non-increasing order, denoted as 
\[
\lambda=\lambda_1\lambda_2\cdots\vdash n,
\quad\text{where $\lambda_1\geq \lambda_2\geq \cdots \geq 1$}.
\]

A \emph{symmetric function} over a commutative ring~$R$ with identity is a formal power series $f(x_1,x_2,\dots)$ of bounded degree such that
\[
f(x_1,x_2,\dots)=f(x_{\sigma(1)},x_{\sigma(2)},\dots)
\]
for every permutation $\sigma$ of the positive integers fixing all but finitely many integers. See \citet[Chapter 7]{Sta99B} for an introduction to symmetric functions. A symmetric function is \emph{homogeneous of degree $d$} if every monomial appearing in it with a nonzero coefficient has total degree~$d$. The \emph{elementary symmetric function} indexed by a partition $\lambda$ is 
\[
e_\lambda=e_{\lambda_1}e_{\lambda_2}\dotsm,
\]
where~$e_k$ is the sum of all products of $k$ distinct variables~$x_i$. Elementary symmetric functions form a basis of the algebra of symmetric functions. A symmetric function is \emph{$e$-positive} if all its $e$-coefficients are nonnegative. We will use the extended symmetric function $e_I$ with composition index $I$, which is defined as 
\[
e_I=e_{\lambda(I)},
\]
where $\lambda(I)$ is the underlying partition whose parts are those of $I$. For a symmetric function $f$, any expansion
\[
f=\sum_{I\vDash n}c_Ie_I
\]
in terms of $e_I$ is called an \emph{$e_I$-expansion} of $f$. An $e_I$-expansion is said to be \emph{positive} if $c_I\ge 0$ for all compositions $I$. Any symmetric function~$f$ having a positive $e_I$-expansion is $e$-positive.

\citet[Definition~2.1]{Sta95} introduced the chromatic symmetric function of a graph $G$, defining it as the formal power series
\[
X_G
=\sum_\kappa 
x_{\kappa(v_1)}\dotsm
x_{\kappa(v_n)},
\]
where $\kappa$ runs over proper colorings of~$G$, and $v_1,\dots,v_n$ are the vertices of $G$. For instance, the complete graph $K_n$ has chromatic symmetric function 
\[
X_{K_n}=n!e_n.
\]

We introduce the functions $w_I'$ and $w_I$ defined for compositions $I$ by
\[
w_I'=(i_2-1)(i_3-1)\dotsm
\quad\text{and}\quad
w_I=i_1 w_I'.
\]
Let $P_n$ denote the $n$-vertex path. \citet[Table~1]{SW16} obtained the following formula.
\begin{lemma}[\citeauthor{SW16}]\label{lem:path}
$X_{P_n}=\sum_{I\vDash n}w_I e_I$ for any $n\ge 1$.
\end{lemma}

For instance, 
\[
X_{P_3}
=
w_3
e_3
+
w_{21}
e_{21}
+
w_{12}
e_{12}
+
w_{111}
e_{111}
=
3e_3
+
e_{12}.
\]

\citet[Lemma 4.4]{Zhe22} presented the following formula for the chromatic symmetric functions of $3$-legged spiders.

\begin{lemma}[\citeauthor{Zhe22}]\label{lem:Sabc}
Let $n\ge 4$. For any partition $(a,b,c)\vdash n-1$, 
\[
X_{S(a,b,c)}
=
X_{P_n}
+\sum_{i=1}^c
(X_{P_i}
X_{P_{n-i}}
-X_{P_{b+i}}
X_{P_{n-b-i}}).
\]
\end{lemma}

\section[Confirming the e-positivity of spiders S(4m+2,2m,1)]{Confirming the \texorpdfstring{$e$}{e}-positivity of spiders \texorpdfstring{$S(4m+2,\,2m,\,1)$}{S(4m+2,2m,1)}}\label{sec:proof}

Throughout this section, we let $m\ge 1$, $n=6m+4$, $G=S(4m+2,\,2m,\,1)$, and 
\[
\mathcal C_d
=\{I\vDash d
\colon i_1, i_2,\dots \ge 2
\}\quad\text{for any $d$.}
\]
Our goal is to find a positive $e_I$-expansion of $X_G$. First, we derive an $e_I$-expansion for $X_G$.

\begin{lemma}\label{lem:X0}
For any $m\ge 1$,
\begin{align*}
X_{S(4m+2,\,2m,\,1)}
&=
e_1^2\sum_{
I\in\mathcal C_{n-2},\
\Theta_I (4m+2)\ge 1}
w_{1I}
e_I
+e_1
X_1
+X_0,
\end{align*}
where
\begin{align}
\notag
X_1
&=\sum_{\substack{
I\in \mathcal C_{n-1}\\
\Theta_I (4m+2)\ge 2}}
w_I
e_I
+
\sum_{\substack{
I\in \mathcal C_{n-1}\\
\Theta_I(4m+2)\ge 1}}
w_{1I}
e_I
+
\sum_{\substack{
J\in\mathcal C_{4m+2}\\ 
K\in\mathcal C_{2m+1}}}
(k_1j_1-2j_1+1)
w_J'
w_K'
e_{J\!K},
\quad\text{and}\\
\creflabel[def]{def:X0}
X_0
&=\sum_{I\in \mathcal C_n} 
w_I
e_I
-
\sum_{(P, Q)\in\mathcal C_{4m+3}\times\mathcal C_{2m+1}}
w_P
w_Q
e_{PQ}.
\end{align}
As a consequence, the spider $S(4m+2,\,2m,\,1)$ is $e$-positive if and only if $X_0$ is $e$-positive.
\end{lemma}
\begin{proof}
Let $G=S(4m+2,\,2m,\,1)$. By \cref{lem:Sabc,lem:path}, 
\begin{align*}
X_G
&=X_{P_n}+X_{P_1}X_{P_{n-1}}-X_{P_{2m+1}}X_{P_{4m+3}}
=\sum_{I\vDash n}
w_I
e_I
+\sum_{I\vDash n-1}
w_I
e_{1I}
-\sum_{J\vDash 4m+3,\
K\vDash 2m+1}
w_J
w_K
e_{J\!K}.
\end{align*}
Arranging the terms according to the number of parts $1$ in compositions, we obtain 
\[
X_G=X_0
+e_1
Y_1
+e_1^2
Y_2,
\]
where
\begin{align}
\label{def:Y1}
Y_1&=
\sum_{
I\in\mathcal C_{n-1}}
(w_{1I}+w_I)e_I
-
\sum_{
J\in\mathcal C_{4m+2},\, 
K\in\mathcal C_{2m+1}}
w_{1J}
w_K 
e_{J\!K}
-
\sum_{
J\in\mathcal C_{4m+3},\,
K\in\mathcal C_{2m}}
w_J
w_{1K}
e_{J\!K}, 
\quad\text{and}\\
\notag
Y_2&=
\sum_{
I\in\mathcal C_{n-2}}
w_{1I}
e_I
-
\sum_{
J\in\mathcal C_{4m+2},\,
K\in\mathcal C_{2m}}
w_{1J}
w_{1K}
e_{J\!K}
=
\sum_{
I\in\mathcal C_{n-2},\,
\Theta_I (4m+2)\ge 1}
w_{1I}
e_I. 
\end{align}
We then regroup the terms in \cref{def:Y1} according to the value of $\Theta_I(4m+2)$ as
\begin{equation}\label{pf:Y1}
Y_1=\sum_{
I\in \mathcal C_{n-1},\
\Theta_I (4m+2)\ge 2}
(w_{1I}+w_I)
e_I
+Z_1
+Z_0,
\end{equation}
where
\begin{align}
\label{def:Z_1}
Z_1
&=\sum_{\substack{
I\in \mathcal C_{n-1}\\
\Theta_I(4m+2)=1}}
(w_{1I}+w_I)
e_I
-\sum_{\substack{
J\in\mathcal C_{4m+3}\\
K\in\mathcal C_{2m}}}
w_J
w_{1K}
e_{J\!K}
=\sum_{\substack{
I\in \mathcal C_{n-1}\\
\Theta_I(4m+2)=1}}
w_{1I}
e_I,
\quad\text{and}\\
\label{def:Z_0}
Z_0
&=
\sum_{\substack{
I\in \mathcal C_{n-1}\\
\Theta_I(4m+2)=0}}
(w_{1I}+w_I)e_I
-
\sum_{\substack{
J\in\mathcal C_{4m+2}\\ 
K\in\mathcal C_{2m+1}}}
w_{1J}
w_K 
e_{J\!K}
=\sum_{\substack{
J\in\mathcal C_{4m+2}\\ 
K\in\mathcal C_{2m+1}}}
(w_{1J\!K}
+w_{J\!K}
-w_{1J}
w_K)
e_{J\!K},
\end{align}
in which the last coefficient can be simplified as 
\[
w_{1J\!K}
+w_{J\!K}
-w_{1J}
w_K
=
(k_1j_1-2j_1+1)w_J' w_K'.
\]
Substituting \cref{def:Z_1,def:Z_0} into \cref{pf:Y1}, we can recast $Y_1$ as the desired symmetric function~$X_1$. Since $j_1,k_1\ge 2$, we have
\[
k_1j_1-2j_1+1\ge 1.
\]
Thus $G$ is $e$-positive if and only if $X_0$ is $e$-positive. 
\end{proof}

Now for any term $w_P w_Q e_{PQ}$ in $X_0$, it is natural to choose 
\[
w_{PQ}e_{PQ} 
\quad\text{and}\quad
w_{QP}e_{QP}
\]
from the positive sum $\sum_{I\in \mathcal C_n}w_I e_I$ and show that 
\[
w_{PQ}
+w_{QP}
-
w_P
w_Q
\ge 0.
\]
But the composition $PQ$ (or $QP$, or both) might be used more than once. In \cref{lem:Y}, we handle the pairs $(P, Q)$ such that each of the compositions $PQ$ and $QP$ is used only once in the positive sum, and characterize the remaining part in $X_0$.

\begin{lemma}\label{lem:Y}
For any $m\ge 1$, 
\[
X_0
=\sum_{
\substack{
(P, Q)\in\mathcal C_{4m+3}\times\mathcal C_{2m+1}\\
\Theta_P(2m+1)\ge 2}}
(p_1q_1-p_1-q_1)
w_P'
w_Q'
e_{PQ}
+Y,
\]
where
\begin{equation}\creflabel[def]{def:Y}
Y
=
\sum_{
(J,K,L)\in \mathbf T}
b(J,K,L)
e_{J\!K\!L}
+
\sum_{I\in \mathcal A}
w_I
e_I,
\end{equation}
in which
\begin{align}
\label{def:A}
\mathcal A
&=\{I\in\mathcal C_n\colon
\Theta_I (2m+1)\Theta_I (4m+3)\ne 0\},\\
\label{def:b(JKL)}
b(J,K,L)
&=w_{J\!K\!L}
+w_{K\!J\!L}
+w_{K\!LJ}
-w_{J\!K} w_L
-w_{K\!J} w_L,
\quad\text{and}\\
\label{def:T}
\mathbf T
&=\mathcal C_{2m+2}\times
\mathcal C_{2m+1}\times
\mathcal C_{2m+1}.
\end{align}
As a consequence, $X_0$ is $e$-positive if $Y$\! is $e$-positive.
\end{lemma}
\begin{proof}
In view of \cref{def:X0}, we separate the negative sum according to the sets
\begin{align*}
\mathbf B_1
&=\{(P,Q)\in \mathcal C_{4m+3}
\times
\mathcal C_{2m+1}\colon
\Theta_P(2m+1)\le 1\}
\quad\text{and}\\
\mathbf B_2
&=\{(P,Q)\in \mathcal C_{4m+3}
\times
\mathcal C_{2m+1}\colon
\Theta_P(2m+1)\ge 2\}.
\end{align*}
Treating $(P,Q)\in \mathbf B_2$, we consider the positive terms with $I\in\mathcal B_2\sqcup\mathcal B_2'$, where
\begin{align}
\label{def:B_2}
\mathcal B_2
&=\{PQ\colon
(P,Q)\in \mathbf B_2\}
=\{I\in\mathcal C_n\colon
\Theta_I (2m+1)\ge 2,\
\Theta_I (4m+2)=1\}
\quad\text{and}\\
\label{def:B_2'}
\mathcal B_2'
&=\{QP\colon
(P,Q)\in \mathbf B_2\}
=\{I\in\mathcal C_n\colon
\Theta_I (2m+1)=0,\
\Theta_I (4m+2)\ge 2
\}.
\end{align}
Then 
\begin{equation}\label{pf:difference2}
\sum_{I\in\mathcal B_2\cup\mathcal B_2'}
w_I
e_I
-
\sum_{
(P,Q)\in \mathbf B_2}
w_P
w_Q
e_{PQ}
=\sum_{
(P,Q)\in \mathbf B_2}
(w_{PQ}
+w_{QP}
-w_P w_Q)
e_{PQ},
\end{equation}
in which the coefficient simplifies to $(p_1q_1-p_1-q_1)w_P'w_Q'$. Next we deal with $\mathcal B_1$. By definition, 
\[
\mathbf B_1
=
\bigsqcup\limits_{(J,K,L)\in
\mathbf T}
\{(J\!K,L),\,
(K\!J,L)\}.
\]
For $(P,Q)\in \mathbf B_1$, we consider the positive terms with $I\in\mathcal B_1\sqcup\mathcal B_1'$, where
\begin{align}
\label{def:B_1}
\mathcal B_1
&=\bigsqcup\limits_{(J,K,L)
\in\mathbf T}
\{J\!K\!L,\,K\!J\!L\}
=\{I\in\mathcal C_n\colon
\Theta_I (2m+1)\le 1,\
\Theta_I (4m+2)=1\}\quad\text{and}\\
\label{def:B_1'}
\mathcal B_1'
&=\bigsqcup\limits_{(J,K,L)
\in\mathbf T}
\{K\!LJ\}
=\{I\in\mathcal C_n\colon
\Theta_I (2m+1)
=\Theta_I (4m+2)
=0\}.
\end{align}
Then
\begin{equation}\label{pf:difference1}
\sum_{I\in \mathcal B_1\sqcup \mathcal B_1'}
w_I
e_I
-
\sum_{
(P,Q)\in \mathbf B_1}
w_P
w_Q
e_{PQ}
=
\sum_{
(J,K,L)\in \mathbf T}
b(J,K,L)
e_{J\!K\!L}.
\end{equation}
In view of \cref{def:B_2,def:B_2',def:B_1,def:B_1'}, the sets $\mathcal B_1$, $\mathcal B_1'$, $\mathcal B_2$ and $\mathcal B_2'$ are pairwise disjoint, and 
\[
\mathcal C_n\backslash \mathcal B_1
\backslash  \mathcal B_1'
\backslash \mathcal B_2
\backslash \mathcal B_2'
=\mathcal A.
\]
Adding up \cref{pf:difference1,pf:difference2}, we obtain the desired formula. Since $p_1,q_1\ge 2$, we find
\[
p_1q_1-p_1-q_1
=(p_1-1)(q_1-1)-1
\ge 0,
\]
and hence complete the proof.
\end{proof}

We remark that a partition $\lambda$ may appear in both parts
\[
\sum_{
\substack{
(P, Q)\in\mathcal C_{4m+3}\times\mathcal C_{2m+1}\\
\Theta_P(2m+1)\ge 2}}
(p_1q_1-p_1-q_1)
w_P'
w_Q'
e_{PQ}
\quad\text{and}\quad
Y
=
\sum_{
(J,K,L)\in \mathbf T}
b(J,K,L)
e_{J\!K\!L}
+
\sum_{I\in \mathcal A}
w_I
e_I
\]
of $X_0$ in \cref{lem:Y}. For example, for $m=10$, the partition $18^33^22^2\vDash 64$ appears in the first sum as the composition $PQ$ with
\[
(P, Q)=\brk1{(18, 18, 3, 2, 2), \,(18, 3)},
\]
and appears in $Y$ as the composition $JKL$ with
\[
(J, K, L)=\brk1{(18, 2, 2), \,(18, 3), \,(18, 3)}.
\]

To find a positive $e_I$-expansion for $Y$, we give a positive $e_I$-expansion for a part of $Y$, and characterize the remaining part, say, $W$. We will give a positive $e_I$-expansion for $W$ later.

Using the symmetry between $K$ and $L$ as compositions in the triple 
\[
(J, K, L)\in\mathbf T
=\mathcal C_{2m+2}\times
\mathcal C_{2m+1}\times
\mathcal C_{2m+1}, 
\] 
we use the lexicographic ordering for compositions: we write $K>L$ if there exists 
\[
1\le s\le \min\brk1{\ell(K),\ell(L)}
\]
such that $k_u=l_u$ for $u\le s-1$, and $k_s>l_s$. For any triple $(j, k, l)$ of numbers, let 
\begin{equation}\creflabel[def]{def:f}
f(j,k,l)
=2j k l
-3j k
-3j l
-2k l
+2j
+2k
+2l.
\end{equation}
For compositions $K$ and $L$ with $K\ge L$, define 
\begin{align*}
g_{K,L}(j,k,l)
&=\begin{dcases*}
f(j,k,l),
& if $K>L$,\\
f(j,k,l)/2,
& if $K=L$.
\end{dcases*}
\end{align*}
For notational simplicity, whenever the pair $(K,L)$ is clear from the context, we write
\[
g(j,k,l)=g_{K,L}(j,k,l).
\]

\begin{lemma}\label{lem:T1234}
The symmetric function $Y$ can be recast as
\begin{equation}\label{YW}
Y
=\sum_{
(J,K,L)\in \mathbf T,\
K\ge L,\
g(j_1,\,k_1,\,l_1)\ge 0}
g(j_1,\,k_1,\,l_1)
w_J'
w_K'
w_L'
e_{J\!K\!L}
+W,
\end{equation}
where $\mathbf T$ is defined by \cref{def:T} and
\[
W=
\sum_{I\in \mathcal A}
w_I
e_I
+
\sum_{
(J,K,L)\in \mathbf T'}
g(j_1,k_1,l_1)
w_J'
w_K'
w_L'
e_{J\!K\!L},
\]
in which $\mathcal A$ is defined by \cref{def:A}, $\mathbf T'=\bigsqcup_{i=1}^4 \mathbf T_i$, and
\begin{align*}
\mathbf T_1
&=\{(J,K,L)\in \mathbf T\colon
l_1=3,\,
k_1=3,\,
j_1=2,\,
K\ge L\},\\
\mathbf T_2
&=\{(J,K,L)\in \mathbf T\colon
l_1=2,\,
f(j_1,k_1,2)<0,\,
\text{$k_1$ is even, $k_1\ge 4$}
\},\\
\mathbf T_3
&=\{(J,K,L)\in \mathbf T\colon
l_1=2,\,
f(j_1,k_1,2)<0,\,
\text{$k_1$ is odd}
\},
\quad\text{and}\\
\mathbf T_4
&=\{(J,K,L)\in \mathbf T\colon
l_1=2,\,
k_1=2,\,
K\ge L\}.
\end{align*}
As a consequence, $Y$\,is $e$-positive if $W$\,is $e$-positive.
\end{lemma}
\begin{proof}
We start from the expression for $Y$ in \cref{def:Y} by simplifying 
\[
\sum_{
(J,K,L)\in \mathbf T}
b(J,K,L)
e_{J\!K\!L}
\]
for $K=L$ and $K\neq L$, respectively. Note that
\[
\sum_{
(J,K,L)\in \mathbf T,\,K\ne L}
b(J,K,L)
e_{J\!K\!L}
=
\sum_{
(J,K,L)\in \mathbf T,\,K>L}
\brk1{b(J,K,L)+b(J,L,K)}
e_{J\!K\!L},
\]
in which
\[
b(J,K,L)
+b(J,L,K)
=
f(j_1,\,k_1,\,l_1)
w_J'
w_K'
w_L',
\]
where $f$ is defined by \cref{def:f}. Then the terms for $K=L$ can be recast as
\[
\sum_{
(J,K,L)\in \mathbf T,\,K=L}
b(J,K,L)
e_{J\!K\!L}
=\frac{1}{2}
\sum_{
(J,K,L)\in \mathbf T,\,
K=L}
f(j_1,\,k_1,\,l_1)
w_J'
w_K'
w_L'
e_{J\!K\!L}.
\]
It follows that 
\[
Y
=\sum_{
(J,K,L)\in \mathbf T,\,K\ge L}
g(j_1,\,k_1,\,l_1)
w_J'
w_K'
w_L'
e_{J\!K\!L}
+
\sum_{I\in \mathcal A}
w_I
e_I.
\]
Comparing it with the desired \cref{YW}, it suffices to show 
\begin{equation}\creflabel[def]{def:T'}
\mathbf T'
=\{
(J,K,L)\in\mathbf T\colon
K\ge L,\,
f(j_1,\,k_1,\,l_1)<0
\}.
\end{equation}

We claim that
\begin{equation*}
f(j,k,l)\ge 0
\quad\text{unless $(j,k,l)=(2,3,3)$ or $\min\{k,\,l\}=2$}.
\end{equation*}
In fact, it is routine to verify that
\[
f(j,k,l)=0
\quad\text{for $(j,k,l)\in
\{(3,3,3),\,(2,3,4),\,(2,4,3)\}$}.
\]
We observe that the function $f$ is increasing in~$k$, since
\[
\frac{\partial f(j,k,l)}
{\partial k}
=(j-1)(2l-3)-1
\ge 0.
\]
It is then increasing in $l$ by symmetry. On the other hand, since
\[
\frac{\partial f(j,k,l)}
{\partial j}
=\frac{(2k-3)(2l-3)-5}{2},
\] 
the function $f$ is increasing in $j$ when $k+l\ge 6$. This proves the claim, and \eqref{def:T'} follows.
\end{proof}

It remains to find a positive $e_I$-expansion for $W$. We will find injections 
$\varphi_i\colon \mathbf T_i\to \mathcal A$ for $1\le i\le 3$ such that
\[
\sum_{I\in \varphi_i(\mathbf T_i)}
w_I
e_I
+
\sum_{(J,K,L)\in \mathbf T_i}
g(j_1,\,k_1,\,l_1)
w_J'
w_K'
w_L'
e_{J\!K\!L}
\]
is $e$-positive and the sets
$\varphi_1(\mathbf T_1)$, 
$\varphi_2(\mathbf T_2)$ and
$\varphi_3(\mathbf T_3)$
are pairwise disjoint, then show that
\[
\sum_{I\in \mathcal A\backslash \varphi_1(\mathbf T_1)\backslash \varphi_2(\mathbf T_2)\backslash \varphi_3(\mathbf T_3)}
w_I
e_I
+
\sum_{(J,K,L)\in \mathbf T_4}
g(j_1,\,k_1,\,l_1)
w_J'
w_K'
w_L'
e_{J\!K\!L}
\]
is $e$-positive. 

As will be seen, for each triple $(J,K,L)\in\mathbf T_i$,
\begin{enumerate}[topsep=2pt]
\item\label[itm]{item1}
$e_{J\!K\!L}=e_{\varphi_i(J,K,L)}$, 
\item\label[itm]{item2}
$c_i(J,K,L)\ge 0$, where
\begin{equation}\creflabel[def]{def:ci}
c_i(J,K,L)
=\frac{w_{\varphi_i(J,K,L)}}{w_J' w_K' w_L'}
+f(j_1,k_1,l_1).
\end{equation}
\end{enumerate}
Since $f(j,k,l)\le g(j,k,l)<0$, the condition \cref{item2} above would guarantee the $e$-positivity of the following part in $W$:
\begin{equation}\label{def:c(JKL)}
w_{\varphi_i(J,K,L)}
e_{\varphi_i(J,K,L)}
+g(j_1,k_1,l_1)
w_J' w_K' w_L'
e_{J\!K\!L}.
\end{equation}
 
Explicitly, we will verify the conditions \cref{item1,item2} above, and give a superset~$\mathcal S_i\supseteq\varphi_i(\mathbf T_i)$ , see \cref{lem:varphi1,lem:varphi2,lem:varphi3}. For the case $i=4$, we construct the map $\varphi_4$ in two pieces and provide supersets $\mathcal S_{41}$ and $\mathcal S_{42}$ for the corresponding images; see \cref{lem:varphi4}. Finally, we shall show in \cref{thm:main} that these supersets are pairwise disjoint, which completes our proof for \cref{conj:epos:4m+2.2m.1}.

Let us start by defining $\varphi_1$. For any 
\[
(J,K,L)\in\mathbf T_1=\{(J,K,L)\in \mathbf T\colon
l_1=3,\,
k_1=3,\,
j_1=2,\,
K\ge L\},
\] 
define $\varphi_1\colon
\mathbf T_1
\to \mathcal A$
 by 
\begin{equation}\creflabel[def]{def:varphi1}
\varphi_1(J,K,L)=P\!K\!R,
\end{equation}
where 
\[
R=j_1\cdot l_2\dotsm l_{-1}
\]
is obtained from $L$ by replacing $l_1=3$ with $j_1=2$, and
\[
P=\min\{j_{-1},\,l_1\}\cdot
j_2\dotsm j_{-2}\cdot
\max\{j_{-1},\,l_1\}.
\]
In other words, the composition 
\[
\varphi_1(J,K,L)=\varphi_1(2j_2\dotsm j_{-1},\,K,\,3l_2\dotsm l_{-1})
\]
is obtained by exchanging $j_{-1}$ and $l_1=3$ if $j_{-1}=2$, and by exchanging $j_1=2$ and $l_1=3$ otherwise, see \cref{fig:1}. We will verify that $\varphi_1(\mathbf T_1)\subseteq \mathcal A$ in
\cref{lem:varphi1}. 
\begin{figure}[H]
\begin{tikzpicture}
\def\xA{0}    
\def\xB{2.2}
\def\xC{2.9}  

\node at (\xA,0)   {$2$};
\node at (\xA+.5,0)   {$j_2$};
\node at (\xA+1,0) {$\dots$};
\node at (\xA+1.5,0) {$j_{-1}$};
\node at (\xB,0)   {$K$};
\node at (\xC,0)   {$3$};
\node at (\xC+.5,0) {$l_2$};
\node at (\xC+1,0) {$\dots$};
\node at (\xC+1.5,0) {$l_{-1}$};

\draw[-] (\xA,.6) -- (\xC,.6);
\draw[->] (\xC,.6) -- (\xC,.25);
\draw[->] (\xA,.6) -- (\xA,.25); 
\draw[-] (\xA+1.5,-.6) -- (\xC,-.6);
\draw[<-] (\xC,-.25) -- (\xC,-.6);
\draw[<-] (\xA+1.5,-.25) -- (\xA+1.5,-.6);
\node at (1.45,.8)   {if $j_{-1}\neq 2$};
\node at (2.2,-.8)   {if $j_{-1}=2$};
\end{tikzpicture}
\caption{The composition $\varphi_1(J,K,L)$ is obtained by exchanging~$j_{-1}$ and $l_1=3$ if $j_{-1}=2$, and by exchanging $j_1=2$ and $l_1=3$ otherwise.}
\label{fig:1}
\end{figure}
It is direct to see that 
\[
e_{\varphi_1(J,K,L)}=e_{J\!K\!L}.
\]
On the other hand, since $f(2,3,3)=-2$, 
\[
c_1(J,K,L)
=\begin{dcases*}
2(3-1)(3-1)-2=6,
& if $j_{-1}=2$\\
3(2-1)(3-1)-2=4,
& if $j_{-1}\ge 3$
\end{dcases*}
\]
is positive. For example,
$\varphi_1(26,\,34,\,322)=3634222$
and
$\varphi_1(242,\,34,\,322)=24334222$.

\begin{lemma}[Properties of $\varphi_1$]\label{lem:varphi1}
We have the following.
\begin{enumerate}
\item\label[itm]{item11}
For any $(J,K,L)\in\mathbf T_1$, we have
\[
\varphi_1(J,K,L)\in\mathcal A\cap\mathcal S_1,
\]
where 
\[
\mathcal S_1
=\{PQR\in \mathcal C_n\colon
(\abs{P}, 
\abs{Q}, 
\abs{R})
=(2m+3,\,
2m+1,\, 
2m)\},
\]
in which the composition triple $(P,Q,R)$ is exactly $(P,K,R)$ that is defined in \cref{def:varphi1}.
\item\label[itm]{item12}
The map $\varphi_1$ is injective.
\end{enumerate}
\end{lemma}
\begin{proof}
Let $I=\varphi_1(J,K,L)=P\!K\!R$. First of all, we check $I\in \mathcal A$. Recall that 
\[
\mathcal A
=\{I\in\mathcal C_n\colon
\Theta_I (2m+1)\Theta_I (4m+3)\ne 0\}.
\]
In fact, we have
\[
\Theta_I (2m+1)=2,
\] 
since $\abs{P}=\abs{J}-j_1+l_1=2m+3$ and since $p_{-1}\ge l_1=3$ by definition. On the other hand, we have 
\[
\Theta_I (4m+3)=1,
\] 
since $\abs{P}+\abs{K}=4m+4$ and since $k_{-1}\ge 2$. This proves $I\in \mathcal A$. It follows that 
\[
\abs{R}=2m
\quad\text{and}\quad
I\in\mathcal S_1.
\]
This proves~\cref{item11}.

Next, we shall show that $\varphi_1$ is injective. Suppose that 
\[
(J',K',L')\in\mathbf T_1
\quad\text{and}\quad
P\!K\!R
=\varphi_1(J,K,L)
=\varphi_1(J',K',L')
=P'\!K'\!R'.
\]
Let $I=P\!K\!R$. Since each of $P$ and $P'$ is the prefix of $I$ of size $2m+3$, we deduce that 
\[
P=P'.
\] 
Similarly, we deduce that 
\[
K=K' 
\quad\text{and}\quad
R=R'.
\]
Since $L$ (resp., $L'$) can be obtained from $R$ (resp., $R'$) by replacing the first part with $3$, we find 
\[
L=L'.
\]
Since $P=P'$, we find 
\[
\ell(J)=\ell(J'), 
\quad\text{and}\quad
j_s=j_s'\ \text{for $2\le s\le \ell(J)-1$};
\]
moreover,
\[
\min\{j_{-1},\,3\}
=\min\{j'_{-1},\,3\}
\quad\text{and}\quad
\max\{j_{-1},\,3\}
=\max\{j'_{-1},\,3\}.
\]
It follows that 
\[
j_{-1}=j_{-1}'.
\]
Since $j_1=2=j_1'$, we find 
\[
J=J'.
\]
This proves the injectivity of $\varphi_1$ and completes the proof.
\end{proof}

Secondly, let us define $\varphi_2$. For any 
\[
(J,K,L)\in\mathbf T_2
=\{(J,K,L)\in \mathbf T\colon
l_1=2,\,
f(j_1,k_1,2)<0,\,
\text{$k_1$ is even, $k_1\ge 4$}
\},
\]
define $\varphi_2\colon\mathbf T_2\to \mathcal A$ by
\begin{equation}\creflabel[def]{def:varphi2}
\varphi_2(J,K,L)=JQR,
\quad\text{where 
$Q=2k_2\dotsm k_{-1}k_1$ and 
$R=l_2\dotsm l_{-1}$}.
\end{equation}
In other words, $\varphi_2(J,K,L)$ is obtained from $J\!K\!L$ by exchanging $k_1$ and $l_1$,
 see \cref{fig:2}. We will verify that $\varphi_2(\mathbf T_2)\subseteq \mathcal A$ in
\cref{lem:varphi2}. 
\begin{figure}[H]
\begin{tikzpicture}
\def\xA{0}    
\def\xB{.7}
\def\xC{2.9}  
    
\node at (\xA,0)   {$J$};
\node at (\xB,0)   {$k_1$};
\node at (\xB+.5,0) {$k_2$};
\node at (\xB+1,0) {$\dots$};
\node at (\xB+1.5,0) {$k_{-1}$};
\node at (\xC,0)   {$2$};
\node at (\xC+.5,0) {$l_2$};
\node at (\xC+1,0) {$\dots$};
\node at (\xC+1.5,0) {$l_{-1}$};

\draw[-] (\xB,.6) -- (\xC,.6);
\draw[->] (\xC,.6) -- (\xC,.25);
\draw[->] (\xB,.6) -- (\xB,.25); 
\end{tikzpicture}
\caption{The composition $\varphi_2(J,K,L)$ is obtained by exchanging~$k_1$ and $l_1=2$.}
\label{fig:2}
\end{figure} 
It is direct to see that 
\[
e_{\varphi_2(J,K,L)}=e_{J\!K\!L}.
\]
On the other hand, 
\begin{equation}\creflabel[ineq]{pf:F'>=2}
c_2(J,K,L)
=j_1(k_1-1)(2-1)
+f(j_1,k_1,2)
=(j_1-1)(2k_1-5)-1
\ge 2.
\end{equation}
For example, $\varphi_2(8,43,232)=823432$.

\begin{lemma}[Properties of $\varphi_2$]\label{lem:varphi2}
We have the following.
\begin{enumerate}
\item
For any $(J,K,L)\in\mathbf T_2$, we have
\[
\varphi_2(J,K,L)\in\mathcal A\cap\mathcal S_2,
\]
where
\[
\mathcal S_2
=\{PQR\in \mathcal C_n\colon
(\abs{P},\abs{Q},\abs{R})
=(2m+2,\,
2m+3,\,
2m-1),\,
q_1=2,\,
\text{$q_{-1}$ is even and $q_{-1}\ge 4$}
\},
\]
in which the composition triple $(P,Q,R)$ is exactly $(J,Q,R)$ that is defined in \cref{def:varphi2}.
\item
The map $\varphi_2$ is injective.
\end{enumerate}
\end{lemma}

\begin{proof}
Let $I=\varphi_2(J,K,L)=JQR$. 
Let us check $I\in \mathcal A\cap\mathcal S_2$. In fact, we have
\[
\Theta_I (2m+1)=1,
\] 
since $\abs{J}=2m+2$ and $j_{-1}\ge 2$. On the other hand, we have
\[
\Theta_I (4m+3)=2,
\] 
since $\abs{J}+\abs{Q}=4m+5$ and $q_{-1}=k_1\ge 4$. It follows that 
\[
\abs{Q}=2m+3
\quad\text{and}\quad
\abs{R}=2m-1.
\]
Moreover, $q_1=l_1=2$, $q_{-1}=k_1$ is even and is at least $4$. This verifies 
\[
I\in\mathcal A\cap\mathcal S_2.
\]

Next, we will show that $\varphi_2$ is injective. Suppose that 
\[
(J',K',L')\in\mathbf T_2
\quad\text{and}\quad
JQR
=\varphi_2(J,K,L)
=\varphi_2(J',K',L')
=J'Q'R'.
\]
Let $I=JQR$. Since $I\in\mathcal S_2$, the composition $J$ is the prefix of $I$ of size $2m+2$, and so is $J'$. Thus 
\[
J=J'.
\]
For the same reason, we have 
\[
Q=Q'
\quad\text{and}\quad
R=R'.
\]
By definition, the composition $L$ (resp., $L'$) can be obtained by inserting a part $2$
at the beginning of $R$ (resp., $R'$). Since $R=R'$, we find 
\[
L=L'.
\]
On the other hand, the composition $K$ (resp., $K'$) can be obtained from $Q$ (resp., $Q'$) by removing the first part and then moving the last part to the beginning. Since $Q=Q'$, we derive that 
\[
K=K'.
\] 
This proves the injectivity of $\varphi_2$ and completes the proof.
\end{proof}

Thirdly, let us define $\varphi_3$. For any triple
\[
(J,K,L)\in\mathbf T_3
=\{(J,K,L)\in \mathbf T\colon
l_1=2,\,
f(j_1,k_1,2)<0,\,
\text{$k_1$ is odd}
\},
\]
each of $K$ and $L$ contains at least one odd part since $\abs{K}=\abs{L}=2m+1$. We call a composition \emph{even} if it has only even parts. For any composition $I$, if $I$ has at least one odd part, let $\fo(I)$ (resp., $\lo(I)$) be the index of the first (resp., last) odd part of $I$, and we use $i_{\fo}$ (resp., $i_{\lo}$) as a shorthand for $i_{\fo(I)}$ (resp., $i_{\lo(I)}$), which is the first (resp., last) odd part of $I$. If $I$ has no odd parts, then the expressions $\fo(I)-1$ and $\ell(I)-\lo(I)$, whenever they occur, are both understood to be $\ell(I)$; in particular, they are $0$ for $I=\emptyset$. Define $\varphi_3\colon\mathbf T_3\to \mathcal A$ by
\begin{equation}\creflabel[def]{def:varphi3}
\varphi_3(J,K,L)=JQR,
\end{equation}
where 
\[
Q=k_2\dotsm k_{-1}
\]
is obtained from $K$ by removing the first part $k_1$, and $R$ is obtained by inserting~$k_1$ into $L$ as:
\begin{equation}\creflabel[def]{eq:R}
R=
2l_2\dotsm l_s\cdot
k_1\cdot
l_{s+1}\dotsm l_{-1}, 
\end{equation}
where $s=\min\{\fo(L)-1,\,\ell(K)-\lo(K)+1\}$. In \cref{eq:R}, the value $\fo(L)-1$ is selected to ensure that $k_1$ is the first odd part of $R$, and the value $\ell(K)-\lo(K)+1$ is selected for proving \cref{eq:eQ}, as will be seen. In other words, the composition $\varphi_3(J,K,L)$ is obtained in the following $2$ steps: 
\begin{enumerate}
\item 
Move the part $k_1$ into $L$ between the parts $l_{\fo(L)-1}$ and $l_{\fo}$.
\item 
If the part $l_{\ell(K)-\lo(K)+1}$ lies to the left of the part $l_{\fo}$, then move $k_1$ into $L$ between the parts 
\[
l_{\ell(K)-\lo(K)+1}
\quad\text{and}\quad
l_{\ell(K)-\lo(K)+2}.
\]
\end{enumerate}
See \cref{fig:3}. We will verify that $\varphi_3(\mathbf T_3)\subseteq \mathcal A$ in \cref{lem:varphi3}. 
\begin{figure}[H]
\begin{tikzpicture}[scale=1.5]
\def\xA{0}    
\def\xB{.5}
\def\xC{4.2}  
    
\node at (\xA,0)   {$J$};
\node at (\xB,0)   {$k_1$};
\node at (\xB+.5,0) {$\dots$};
\node at (\xB+1.2,0) {$k_{-x-1}$};
\node at (\xB+2,0) {$k_{-x}$};
\node at (\xB+2.6,0) {$\dots$};
\node at (\xB+3.1,0) {$k_{-1}$};
\node at (\xC,0)   {$2$};
\node at (\xC+.5,0) {$l_2$};
\node at (\xC+1,0) {$\dots$};
\node at (\xC+1.5,0) {$l_y$};
\node at (\xC+2.15,0) {$l_{y+1}$};
\node at (\xC+2.8,0) {$\dots$};
\node at (\xC+3.45,0) {$l_{x+1}$};
\node at (\xC+4.1,0) {$\dots$};
\node at (\xC+4.6,0) {$l_{-1}$};

\draw[-] (\xB,.6) -- (\xC+1.65,.6);
\draw[->] (\xC+1.65,.6) -- (\xC+1.65,.25);
\draw[-] (\xB,.6) -- (\xB,.25); 
\draw[dashed] (\xB+3.45,.5) -- (\xB+3.45,-.5); 
\draw[-] (\xC+1.75,.5) -- (\xC+1.75,-.5); 
\begin{scope}[yshift=-1.5cm]
\def\xA{0}    
\def\xB{.5}
\def\xC{4.2}  
    
\node at (\xA,0)   {$J$};
\node at (\xB,0)   {$k_1$};
\node at (\xB+.5,0) {$\dots$};
\node at (\xB+1.2,0) {$k_{-x-1}$};
\node at (\xB+2,0) {$k_{-x}$};
\node at (\xB+2.6,0) {$\dots$};
\node at (\xB+3.1,0) {$k_{-1}$};
\node at (\xC,0)   {$2$};
\node at (\xC+.5,0) {$l_2$};
\node at (\xC+1,0) {$\dots$};
\node at (\xC+1.65,0) {$l_{x+1}$};
\node at (\xC+2.45,0) {$l_{x+2}$};
\node at (\xC+3.1,0) {$\dots$};
\node at (\xC+3.75,0) {$l_{y+1}$};
\node at (\xC+4.4,0) {$\dots$};
\node at (\xC+4.9,0) {$l_{-1}$};

\draw[-] (\xB,.6) -- (\xC+2.05,.6);
\draw[->] (\xC+2.05,.6) -- (\xC+2.05,.25);
\draw[-] (\xB,.6) -- (\xB,.25); 
\draw[dashed] (\xB+3.45,.5) -- (\xB+3.45,-.5); 
\draw[-] (\xC+3.4,.5) -- (\xC+3.4,-.5); 
\end{scope}
\end{tikzpicture}
\caption{The map $\varphi_3(J,K,L)$.
The upper part corresponds to the case $\fo(L)-1<\ell(K)-\lo(K)+1$, and the lower part corresponds to the case $\fo(L)-1\geq \ell(K)-\lo(K)+1$. Here $x=\ell(K)-\lo(K)$ and $y=\fo(L)-1$.}
\label{fig:3}
\end{figure}
Since $\varphi_3(J,K,L)$ is obtained from $J\!K\!L$ by moving~$k_1$ rightward, we find 
\[
e_{\varphi_3(J,K,L)}=e_{J\!K\!L}.
\]
As in \cref{pf:F'>=2}, we have
\[
c_3(J,K,L)
=j_1(k_1-1)+f(j_1,k_1,2)
=(j_1-1)(2k_1-5)-1\ge 0,
\]
since $j_1\ge 2$ and $k_1\ge 3$ is odd. For example,
\[
\varphi_3(86,3433,2326)=8643323326\quad\text{and}\quad\varphi_3(86,34222,2326)=86422223326.
\]

\begin{lemma}[Properties of $\varphi_3$]\label{lem:varphi3}
We have the following.
\begin{enumerate}
\item
For any $(J,K,L)\in\mathbf T_3$, we have 
\[
\varphi_3(J,K,L)\in\mathcal A\cap\mathcal S_3,
\]
where
\begin{multline*}
\mathcal S_3
=
\{
PQR\in\mathcal C_n\colon
r_1=2,\,
\text{$R$ has at least $2$ odd parts},\\
(\abs{P},
\abs{Q},
\abs{R})
=(2m+2,\,
2m+1-r_{\fo},\,
2m+1+r_{\fo}),\\ 
\text{
either $\fo(R)-1\le \ell(Q)-\lo(Q)$ and $r_{\fo(R)+1}$ is odd, or 
$\fo(R)-1=\ell(Q)-\lo(Q)+1$}\},
\end{multline*}
in which the composition triple $(P,Q,R)$ is exactly $(J,Q,R)$ that is defined in \cref{def:varphi3}.
\item
The map $\varphi_3$ is injective.
\end{enumerate}
\end{lemma}
\begin{proof}
Let $I=\varphi_3(J,K,L)=JQR$. Let us check $I\in \mathcal A$. In fact, we have
\[
\Theta_I (2m+1)=1,
\] 
since $\abs{J}=2m+2$ and $j_{-1}\ge 2$ by definition. On the other hand, assume to the contrary that $\Theta_I(4m+3)=0$. Since
\[
\abs{J}+\abs{Q}=4m+3-k_1,
\]
the composition $R$ has a prefix of size $k_1$.
By the definition of $R$, the prefix ending at the inserted part $k_1$ is
\[
l_1\dotsm l_s k_1,
\quad\text{where $s=\min\{\fo(L)-1,\,\ell(K)-\lo(K)+1\}$}.
\]
Moreover, the parts $l_1,\dots,l_s$ are all even. Hence a prefix of
$R$ can have size $k_1$ only if $\abs{l_1\dotsm l_s}=0$, namely
only if $s=0$. This is impossible, since 
\[
\fo(L)-1\ge 1
\quad\text{and}\quad
\ell(K)-\lo(K)+1\ge 1.
\]
This contradiction proves $I\in\mathcal A$.

Next, we will show that $I\in \mathcal S_3$. Since $l_1=2$, we find 
\[
\fo(L)-1\ge 1
\quad\text{and}\quad
r_1=l_1=2.
\]
Since $\abs{L}=2m+1$ is odd and $k_1$ is odd, the composition $R$ has at least two odd parts. From the definition, it is direct to see that $\abs{Q}$ and $\abs{R}$ are of the desired values, and $\abs{P}=n-\abs{Q}-\abs{R}$ is as desired. By definition, we see that 
\[
\ell(Q)-\lo(Q)=\ell(K)-\lo(K).
\]
To see that $I$ satisfies the last condition in $\mathcal S_3$, we proceed according to the definition of $R$. 
\begin{enumerate}
\item If $\fo(L)-1\le \ell(K)-\lo(K)$, then \cref{eq:R} reduces to
\[
R=
l_1\dotsm l_{\fo(L)-1}\cdot
k_1\cdot
l_{\fo}\dotsm l_{-1}.
\]
It follows that 
\[
\fo(R)-1=\fo(L)-1\le \ell(K)-\lo(K)=\ell(Q)-\lo(Q),
\] 
and that 
\[
r_{\fo(R)+1}=l_{\fo},
\]
which is defined to be odd.
\item 
If $\fo(L)-1\ge \ell(K)-\lo(K)+1$, then \cref{eq:R} reduces to
\[
R=
l_1\dotsm l_{\ell(K)-\lo(K)+1}\cdot
k_1\cdot
l_{\ell(K)-\lo(K)+2}\dotsm l_{-1}.
\]
It follows that 
\[
\fo(R)-1=\ell(K)-\lo(K)+1=\ell(Q)-\lo(Q)+1.
\]
\end{enumerate}
This proves $I\in \mathcal S_3$.

Second, let us show that $\varphi_3$ is injective. Suppose that $(J',K',L')\in\mathbf T_3$ and
\[
JQR
=\varphi_3(J,K,L)
=\varphi_3(J',K',L')
=J'Q'R', 
\]
where $Q$, $R$, $Q'$ and $R'$ are defined by \cref{def:varphi3}. Let $I=JQR$. By definition of $\varphi_3$, the composition~$J$ is the prefix of $I$ of size $2m+2$, and so is $J'$. Thus 
\[
J=J'
\quad\text{and}\quad 
QR=Q'R'.
\]
Assume that $Q\ne Q'$. Then $\abs{Q}\ne\abs{Q'}$. Suppose that $\abs{Q}<\abs{Q'}$ without loss of generality. Then there exists a composition $M\ne \emptyset$ such that
\[
Q'=QM\quad\text{and}\quad R=MR'.
\]
Since $I\in \mathcal S_3$,
\[
2m+1-r'_{\fo}
=\abs{Q'}
>\abs{Q}
=2m+1-r_{\fo}.
\]
It follows that 
\[
r_{\fo}>r_{\fo}'.
\]
If $M$ had no odd part, then the first odd part of $R=MR'$ would be the first odd part of $R'$, of value~$r_{\fo}'$, a contradiction. Thus $M$ has an odd part, and since $M$ is a prefix of $R$, the first odd part of $M$ is the first odd part of $R$, of value $r_{\fo}$. This proves
\[
\abs{Q'}
\ge \abs{Q}+r_{\fo}
=2m+1,
\]
a contradiction. Thus 
\[
Q=Q'.
\]
As a result, we obtain 
\[
R=R'.
\]
In summary, we have
\[
(J,Q,R)=(J',Q',R').
\]
Since $K$ (resp., $K'$) is obtained by inserting the first odd part of $R$ (resp., $R'$) at the beginning of~$Q$ (resp., $Q'$), we find 
\[
K=K'.
\]
Since $L$ (resp., $L'$) is obtained from $R$ (resp., $R'$) by removing the first odd part, we find 
\[
L=L'.
\]
This proves the injectivity of $\varphi_3$ and completes the proof.
\end{proof}

Now let us handle $\mathbf T_4$. Recall that
\[
\mathbf T_4
=\{(J,K,L)\in \mathbf T\colon
k_1=l_1=2,\,
K\ge L\}.
\]
We call a composition \emph{odd} if it has only odd parts. For $(J,K,L)\in\mathbf T_4$, let $U=U(K)$ be the composition obtained from $K$ by moving its longest odd suffix to the beginning; in other words, 
\[
U=K_2 K_1,
\]
where $K_2$ is the longest odd suffix of $K$, and $K_1$ is the prefix of $K$ such that $K=K_1K_2$. Then $u_{-1}$ is even, i.e., 
\begin{equation}\creflabel[ineq]{eq:e-(U)}
\ell(U)-\lo(U)\geq 1.
\end{equation}
Note that $K_2$ is empty if $k_{-1}$ is even; in that case $U=K$. For any $(J,K,L)\in\mathbf T_4$, let $I=JUL$. 

Define $\varphi_4(J,K,L)$ to be the composition obtained in the following $3$ steps:
\begin{enumerate}
\item
Insert a bar between the parts $u_{-1}$ and $l_1$ in $I$. 
\item
Locate all odd parts of $I$ that are nearest to the bar. Since $\abs{K}=\abs{L}=2m+1$, such a part always exists. When such a part is not unique, the number of such parts must be two; in that case, we assign the part lying to the right of the bar to be the desired one.
\item
Move this odd part to the symmetric position on the other side of the bar.
\end{enumerate}
For example, 
\begin{align}
\label{eq:varphi41}
\varphi_{41}(86, 2335, 232222)&=863353222222,\quad\text{and}\\
\label{eq:varphi42}
\varphi_{42}(86, 22333, 222232)&=8633222232232.
\end{align}
Since $u_{-1}$ is even, we can infer that the odd part is not the one adjacent to the bar; thus
\begin{equation}\creflabel[ineq]{ineq:JU(K)L}
\varphi_4(J,K,L)\neq JUL. 
\end{equation}
We will need \cref{ineq:JU(K)L} in proving $\varphi_4(J,K,L)\in \mathcal A$; this is also the reason why we introduce the operation $U$. 

For technical convenience, we give an algebraic definition for $\varphi_4$, which is equivalent to the combinatorial one above. Let 
\[
\mathbf T_4'=\{(J,K,L)\in\mathbf T,\,k_1=l_1=2\}.
\]
Then $\mathbf T_4\subset \mathbf T_4'$, and $\mathbf T_4'$ decomposes as
$\mathbf T_4'=\mathbf T_{41}\sqcup\mathbf T_{42}$, where
\begin{align*}
\mathbf T_{41}
&=\brk[c]1{(J,K,L)\in\mathbf T\colon
k_1=l_1=2,\,
\fo(L)-1\le \ell\brk1{U(K)}-\lo\brk1{U(K)}
},\quad\text{and}\\
\mathbf T_{42}
&=\brk[c]1{(J,K,L)\in\mathbf T\colon
k_1=l_1=2,\,
\fo(L)-1>\ell\brk1{U(K)}-\lo\brk1{U(K)}}.
\end{align*}
For any $(J,K,L)\in\mathbf T_{41}$, define the map $\varphi_{41}\colon\mathbf T_{41}\to \mathcal A$ by 
\begin{equation}\creflabel[def]{def:varphi41}
\varphi_{41}(J,K,L)=JQR,
\end{equation}
which is obtained from the composition $J\cdot U(K)\cdot L$ by moving the first odd part $l_{\fo}$  of $L$ into $U(K)$ such that $\ell(Q)-\lo(Q)=\fo(L)-1$, i.e.,
\[
Q=u_1 \dotsm u_{-\fo(L)} \cdot 
l_{\fo} \cdot 
u_{-\fo(L)+1} \dotsm u_{-1}\quad\text{and}\quad R=2l_2\dotsm l_{\fo(L)-1}l_{\fo(L)+1}\dotsm l_{-1}.
\]
See \cref{eq:varphi41} for example. Note that $l_{\fo}$ exists since $\abs{L}=2m+1$ is odd. It is possible that $R=2$. For any $(J,K,L)\in \mathbf T_{42}$, define $\varphi_{42}\colon\mathbf T_{42}\to \mathcal A$ by 
\begin{equation}\creflabel[def]{def:varphi42}
\varphi_{42}(J,K,L)=JQR,
\end{equation}
which is obtained from the composition $J\cdot U(K)\cdot L$ by moving the last odd part $u_{\lo}$  of $U=U(K)$ into $L$ such that $\fo(R)-1=\ell(U)-\lo(U)$, i.e.,
\[
Q=u_1 \dotsm u_{\lo(U)-1}u_{\lo(U)+1} \dotsm u_{-1}
\quad\text{and}\quad
R=2l_2 \dotsm l_{\ell(U)-\lo(U)} \cdot 
u_{\lo} \cdot 
l_{\ell(U)-\lo(U)+1} \dotsm l_{-1}.
\]
See \cref{eq:varphi42} for example. We will verify that $\varphi_{4i}(\mathbf T_{4i})\subseteq \mathcal A$ in \cref{lem:varphi4}. Now we can define $\varphi_4\colon\mathbf T_4'\to\mathcal A$ by
\[
\varphi_4(J,K,L)
=\begin{dcases*}
\varphi_{41}(J,K,L),
& if $(J,K,L)\in \mathbf T_{41}$,\\
\varphi_{42}(J,K,L),
& if $(J,K,L)\in \mathbf T_{42}$.
\end{dcases*}
\]

As will be seen, we shall show that for each triple $(J,K,L)\in\mathbf T_4$,
\begin{enumerate}
\item\label[itm]{item41}
$e_{\varphi_4(J,K,L)}
=e_{\varphi_4(J,L,K)}
=e_{J\!K\!L}$, and
\item\label[itm]{item42}
$c_4(J,K,L)=0$, where $c_4(J,K,L)$ is defined as
\[
c_4(J,K,L)
=\begin{dcases*}
\dfrac{w_{\varphi_4(J,K,L)}}{w_J' w_K' w_L'}
+
\dfrac{w_{\varphi_4(J,L,K)}}{w_J' w_K' w_L'}
+f(j_1,k_1,l_1),
& if $K>L$,\\[6pt]
\dfrac{w_{\varphi_4(J,K,K)}}{w_J' (w_K')^2}
+g(j_1,k_1,k_1),
& if $K=L$,
\end{dcases*}
\]
cf.~\cref{def:ci}.
\end{enumerate}
The condition \cref{item42} above guarantees the $e$-positivity of the relevant part in $W$. If $K>L$, this part is
\begin{equation*}
w_{\varphi_4(J,K,L)}
e_{\varphi_4(J,K,L)}
+
w_{\varphi_4(J,L,K)}
e_{\varphi_4(J,L,K)}
+
g(j_1,k_1,l_1)
w_J' w_K' w_L'
e_{J\!K\!L}.
\end{equation*}
If $K=L$, the two images coincide and only one positive term is used:
\begin{equation*}
w_{\varphi_4(J,K,K)}
e_{\varphi_4(J,K,K)}
+
g(j_1,k_1,k_1)
w_J' (w_K')^2
e_{J\!K\!K}.
\end{equation*}

Now we are going to prove \cref{item41} and \cref{item42}. In fact, for any $(J,K,L)\in\mathbf T_{41}$, the composition $\varphi_{41}(J,K,L)$ is obtained from the composition~$J\cdot U(K)\cdot L$ by moving the part $l_{\fo}$ leftward. Since $U(K)$ is a rearrangement of $K$, we have 
\[
e_{\varphi_{41}(J,K,L)}=e_{J\!K\!L}.
\]
Similarly, for any $(J',K',L')\in\mathbf T_{42}$, the composition $\varphi_{42}(J',K',L')$ is obtained from the composition~$J'\cdot  U(K')\cdot L'$ by moving the part $u_{\lo}$ rightward. Thus 
\[
e_{\varphi_{42}(J',K',L')}=e_{J'\!K'\!L'}.
\] 
This proves \cref{item41}. On the other hand, since $f(j,2,2)=-2j$, we find
\[
g(j,2,2)=\begin{dcases*}
-2j, & if $K>L$,\\
-j, & if $K=L$.
\end{dcases*}
\]
Moreover, the positive contribution is $2j_1$ when $K>L$ and $j_1$ when $K=L$. Hence
\begin{align*}
c_4(J,K,L)
&=
\begin{dcases*}
j_1(2-1)(2-1)\cdotp 2
-2j_1=0,
& if $K>L$\\
j_1(2-1)(2-1)
-j_1=0,
& if $K=L$
\end{dcases*}
\end{align*}
which is always zero. This proves \cref{item42}.

Now let us provide the aforementioned superset of $\varphi_4(\mathbf T_4)$. We do this by giving a superset for $\varphi_{41}(\mathbf T_{41})$ and for $\varphi_{42}(\mathbf T_{42})$ individually.

\begin{lemma}[Properties of $\varphi_4$]\label{lem:varphi4}
We have the following.
\begin{enumerate}
\item
For any $(J,K,L)\in \mathbf T_{41}$, we have
\[
\varphi_{41}(J,K,L)\in \mathcal A\cap \mathcal S_{41},
\]
where
\begin{multline*}
\mathcal S_{41}
=\{
PQR\in\mathcal C_n\colon 
r_1=2,\,
\fo(R)-1\geq \ell(Q)-\lo(Q),\,
\text{$Q$ has at least two odd parts,}\\
(\abs{P},\abs{Q},\abs{R})
=(2m+2,\,
2m+1+q_{\lo},\,
2m+1-q_{\lo})\},
\end{multline*}
in which the composition triple $(P,Q,R)$ is exactly $(J,Q,R)$ that is defined in \cref{def:varphi41}.
\item\label[itm]{S42}
For any $(J,K,L)\in \mathbf T_{42}$, we have
\[
\varphi_{42}(J,K,L)\in\mathcal A\cap \mathcal S_{42},
\]
where
\begin{multline*}
\mathcal S_{42}
=\{
PQR\in\mathcal C_n\colon
r_1=2,\,
\fo(R)-1 \leq \ell(Q)-\lo(Q),\,
\text{$R$ has at least two odd parts},\\
\text{$r_{\fo(R)+1}$ is even},\,
(\abs{P},\abs{Q},\abs{R})
=(2m+2,\,
2m+1-r_{\fo},\,
2m+1+r_{\fo})
\},
\end{multline*}
in which the composition triple $(P,Q,R)$ is exactly $(J,Q,R)$ that is defined in \cref{def:varphi42}.
\item
Both $\varphi_{41}$ and $\varphi_{42}$ are injective.
\end{enumerate}
\end{lemma}
\begin{proof}
We proceed to show them individually.
\begin{enumerate}
\item
Let $I=\varphi_{41}(J,K,L)=JQR$. We first check $I\in \mathcal A$. In fact, we have
\[
\Theta_I (2m+1)=1,
\] 
since $\abs{J}=2m+2$ and $j_{-1}\ge 2$. Assume to the contrary that $\Theta_I (4m+3)=0$. Since
\[
\abs{J}+\abs{Q}=4m+3+l_{\fo},
\]
the composition $Q$ has a suffix $Q'$ of odd size~$l_{\fo}$. On the other hand, by \cref{eq:e-(U)}, and since~$l_{\fo}$ is the last odd part of $Q$, we find 
\[
\abs{Q'}\ge u_{-1}+l_{\fo}>l_{\fo},
\]
a contradiction. This proves $I\in\mathcal A$.

Next, we will show that $I\in \mathcal S_{41}$. Recall that
\[
\mathbf T_{41}
=\brk[c]1{(J,K,L)\in\mathbf T\colon
k_1=l_1=2,\,
\fo(L)-1\le \ell(U(K))-\lo(U(K))
}.
\]
From the definition, we see that 
\[
r_1=l_1=2
\quad\text{and}\quad
\fo(R)-1\ge \fo(L)-1=\ell(Q)-\lo(Q).
\] 
Since $Q$ is obtained by inserting an odd part into the composition $U(K)$ whose size is odd, it contains at least $2$ odd parts. On the other hand, 
\[
\abs{Q}=\abs{U(K)}+l_{\fo}=2m+1+q_{\lo}
\quad\text{and}\quad
\abs{R}
=n-\abs{J}-\abs{Q}
=2m+1-q_{\lo}.
\]
This proves $I\in\mathcal S_{41}$.
\item
Let $(J,K,L)\in \mathbf T_{42}$ and let $I=\varphi_{42}(J,K,L)=JQR$, where $Q$ and $R$ are defined by \cref{def:varphi42}. We first check $I\in \mathcal A$. In fact, we have
\[
\Theta_I (2m+1)=1
\]
since $\abs{J}=2m+2$ and $j_{-1}\ge 2$. Assume to the contrary that $\Theta_I (4m+3)=0$. Since 
\[
\abs{J}+\abs{Q}=4m+3-u_{\lo},
\]
the composition $R$ has a prefix $R'$ of odd size~$u_{\lo}$. Recall that
\[
\mathbf T_{42}
=\{(J,K,L)\in\mathbf T\colon
k_1=l_1=2,\,
\fo(L)-1>\ell(U(K))-\lo(U(K))
\}.
\]
Then 
\[
\fo(L)-1>\ell(U(K))-\lo(U(K))\ge 1,
\quad\text{and}\quad
\abs{R'}\ge l_1+u_{\lo}>u_{\lo},
\]
a contradiction. This proves $I\in\mathcal A$.

Next, we will show that $I\in \mathcal S_{42}$. From the definition, we see that 
\[
r_1=l_1=2
\quad\text{and}\quad
\ell(Q)-\lo(Q)\ge \ell(U(K))-\lo(U(K))=\fo(R)-1.
\] 
Since $R$ is obtained by inserting an odd part into~$L$ whose size is odd, it contains at least $2$ odd parts. On the other hand,
\[
\abs{Q}=\abs{U(K)}-u_{\lo}=2m+1-r_{\fo}
\quad\text{and}\quad
\abs{R}
=n-\abs{J}-\abs{Q}
=2m+1+r_{\fo}.
\]
From the definition, we see that $r_{\fo(R)+1}=l_{\ell(U(K))-\lo(U(K))+1}$ is even. This proves $I\in\mathcal S_{42}$.
\item
Let us show that $\varphi_{41}$ is injective. Suppose that 
\[
JQR
=\varphi_{41}(J,K,L)
=\varphi_{41}(J',K',L')
=J'Q'R',
\]
where $Q$, $R$, $Q'$ and $R'$ are defined by 
\cref{def:varphi41}.
We need to show $(J,K,L)=(J',K',L')$.

We claim that $(J,Q,R)=(J',Q',R')$. Let $I=JQR=J'Q'R'$. Since $I\in \mathcal S_{41}$, 
\begin{align*}
(\abs{J},\,
\abs{Q},\,
\abs{R})
&=(2m+2,\,
2m+1+q_{\lo},\,
2m+1-q_{\lo}),\quad\text{and}\\
(\abs{J'},\,
\abs{Q'},\,
\abs{R'})
&=(2m+2,\,
2m+1+q_{\lo}',\,
2m+1-q_{\lo}').
\end{align*}
Thus the composition $J$ is the prefix of $I$ of size $2m+2$, and so is $J'$. It follows that 
\[
J=J' 
\quad\text{and}\quad
QR=Q'R'.
\]
Assume that $Q\ne Q'$. Then we can suppose that $\abs{Q}>\abs{Q'}$ without loss of generality. It follows that 
\[
q_{\lo}>q_{\lo}'.
\]
Thus the last odd part of $Q$, which is $q_{\lo}$, is not contained in~$Q'$. Then 
\[
\abs{Q'}\le \abs{Q}-q_{\lo}=2m+1,
\]
a contradiction. Hence 
\[
\abs{Q}=\abs{Q'} 
\quad\text{and}\quad
q_{\lo}=q_{\lo}'.
\]
It follows that 
\[
Q=Q'
\quad\text{and}\quad
R=R'.
\]
This proves the claim. 

Alternatively, the value $q_{\lo}=q_{\lo}'$ can be obtained by the following $2$ steps:
\begin{enumerate}
\item 
Insert a bar between two parts $i_{x+1}$ and $i_{x+2}$ in $I$, where $x$ is defined by
\[
\sum_{j=1}^x i_j\leq 4m+3\quad\text{and}\quad\sum_{j=1}^{x+1} i_j\geq 4m+4.
\]
\item
Locate the odd part to the left of the bar that is nearest to the bar. Then the value of this part is $q_{\lo}$.
\end{enumerate}

Next, we shall show $(K,L)=(K',L')$. Since $U(K)$ (resp., $U(K')$) can be obtained from $Q$ (resp., $Q'$) by removing the last odd part, we find 
\[
U(K)=U(K').
\]
Since $K$ (resp., $K'$) can be obtained from $U(K)$ (resp., $U(K')$) by moving its longest odd prefix to the end, we find 
\[
K=K'.
\]
On the other hand, since
\[
R=2\cdot l_2\dotsm l_{\fo(L)-1}\cdot
l_{\fo(L)+1}\dotsm l_{-1},
\]
in which 
\[
\fo(L)-1=\ell(Q)-\lo(Q),
\]
the composition $L$ (resp., $L'$) is obtained by inserting a part of value $l_{\fo}=q_{\lo}$ (resp., $l_{\fo}'=q_{\lo}'$) into~$R$ (resp., $R'$) between 
\[
r_{\ell(Q)-\lo(Q)}
\quad\text{and}\quad
r_{\ell(Q)-\lo(Q)+1}.
\]
Since $q_{\lo}=q_{\lo}'$, we find 
\[
l_{\fo}=l_{\fo}' 
\quad\text{and}\quad
L=L'.
\]
This proves the injectivity of $\varphi_{41}$.

Now let us show that $\varphi_{42}$ is injective. We can confirm in the same way above, algebraically. Here we give a combinatorial proof as aforementioned. Suppose that $I\in\varphi_{42}(\mathbf T_{42})$. By \cref{S42}, the composition $I$ can be factored as $JQR\in\mathcal S_{42}$. We shall show that the triple $(J,Q,R)$ is determined uniquely by the image $I$. Since
\begin{align*}
(\abs{J},\,
\abs{Q},\,
\abs{R})
&=(2m+2,\,
2m+1-r_{\fo},\,
2m+1+r_{\fo}),
\end{align*}
the composition $J$ is the prefix of $I$ of size $2m+2$. Insert a bar between the parts $i_x$ and $i_{x+1}$ in~$I$, where the index $x$ is defined by
\[
\sum_{j=1}^x i_j\leq 4m+3\quad\text{and}\quad\sum_{j=1}^{x+1} i_j\geq 4m+4.
\]
Since $R$ has size at least $2m+4$, it has a suffix $I_1=i_{x+1}\dotsm i_{-1}$. 
Since $I\in\mathcal A$, we have 
\[
\Theta_I(4m+3)\ne0,
\]
and hence no prefix of $I$ has size $4m+3$. 
By the choice of $x$, it follows that
\[
\sum_{j=1}^x i_j\le 4m+2,
\]
and consequently
\[
\abs{I_1}=n-\sum_{j=1}^x i_j\ge (6m+4)-(4m+2)=2m+2.
\]
Suppose that $R=MI_1$. 
If $M$ is not even, then $M$ contains the first odd part $r_{\fo}$ of $R$, and so
\[
\abs{R}=\abs{M}+\abs{I_1}\geq r_{\fo}+2m+2,
\]
which contradicts $\abs{R}=2m+1+r_{\fo}$. 
Thus $M$ is even.
Then the first odd part $r_{\fo}$ of $R$ is the first odd part of $I_1$, i.e., $r_{\fo}$ is the odd part of $I$ that is nearest to the bar and lies to the right of the bar. Hence the triple $(J,Q,R)$ is unique.

Now let us find the triple $(J,K,L)$. The composition $L$ can be obtained from $R$ by removing the first odd part $r_{\fo}$. Note that 
\[
Q=u_1\dotsm u_{\lo(U(K))-1}\cdot
u_{\lo(U(K))+1}\dotsm u_{-1},
\]
in which 
\[
\ell(U(K))-\lo(U(K))=\fo(R)-1.
\]
Therefore, $U(K)$ is obtained from $Q$ by inserting a part of value $u_{\lo}=r_{\fo}$ into $Q$ to the immediate left of the $\fo(R)-1$th last part. The composition~$K$ is obtained from $U(K)$ by moving its longest odd prefix to the end. Hence we obtain the triple $(J,K,L)$. This proves the injectivity of $\varphi_{42}$. 
\end{enumerate}
This completes the proof.
\end{proof}

Now, we are in a position to finish the whole proof.

\begin{theorem}\label{thm:main}
For any $m\ge 1$, the spider $S(4m+2,\,2m,\,1)$ is $e$-positive.
\end{theorem}
\begin{proof}
By \cref{lem:X0,lem:Y,lem:T1234}, it suffices to prove the $e$-positivity of $W$. By \cref{lem:varphi1,lem:varphi2,lem:varphi3,lem:varphi4}, it suffices to check that the sets 
\[
\mathcal S_1,\quad
\mathcal S_2,\quad 
\mathcal S_3,\quad
\mathcal S_{41},\quad\text{and}\quad
\mathcal S_{42}
\]
are pairwise disjoint. Recall that
\begin{align*}
\mathcal S_1
&=\{PQR\in \mathcal C_n\colon
(\abs{P}, 
\abs{Q}, 
\abs{R})
=(2m+3,\,
2m+1,\, 
2m)\},\\
\mathcal S_2
&=\{PQR\in \mathcal C_n\colon
(\abs{P},\abs{Q},\abs{R})
=(2m+2,\,
2m+3,\,
2m-1),\,
q_1=2,\,
\text{$q_{-1}$ is even and $q_{-1}\ge 4$}\},\\
\mathcal S_3
&=\{
PQR\in\mathcal C_n\colon
r_1=2,\,
\text{$R$ has at least $2$ odd parts},\\
&\quad\quad
(\abs{P},
\abs{Q},
\abs{R})
=(2m+2,\,
2m+1-r_{\fo},\,
2m+1+r_{\fo}),\\ 
&\quad\quad
\text{
either $\fo(R)-1\le \ell(Q)-\lo(Q)$ and $r_{\fo(R)+1}$ is odd, or 
$\fo(R)-1=\ell(Q)-\lo(Q)+1$}\},\\
\mathcal S_{41}
&=\{
PQR\in\mathcal C_n\colon 
r_1=2,\,
\fo(R)-1\geq \ell(Q)-\lo(Q),\,
\text{$Q$ has at least two odd parts,}\\
&\quad\quad
(\abs{P},\abs{Q},\abs{R})
=(2m+2,\,
2m+1+q_{\lo},\,
2m+1-q_{\lo})\},\quad\text{and}\\
\mathcal S_{42}
&=\{
PQR\in\mathcal C_n\colon
r_1=2,\,
\fo(R)-1 \leq \ell(Q)-\lo(Q),\,
\text{$R$ has at least two odd parts},\\
&\quad\quad
\text{$r_{\fo(R)+1}$ is even},\,
(\abs{P},\abs{Q},\abs{R})
=(2m+2,\,
2m+1-r_{\fo},\,
2m+1+r_{\fo})
\}.
\end{align*}
We find 
\[
\mathcal S_1\cap
(\mathcal S_2\cup
\mathcal S_3\cup
\mathcal S_{41}\cup
\mathcal S_{42})
=\emptyset.
\]
Indeed, $\Theta_I(2m+2)=1$ for $I\in\mathcal S_1$, whereas 
\[
\Theta_I(2m+2)=0\quad\text{for $I\in
\mathcal S_2\cup
\mathcal S_3\cup
\mathcal S_{41}\cup
\mathcal S_{42}$}.
\]
Below we shall show that $\mathcal S_2$, $\mathcal S_3$, $\mathcal S_{41}$ and $\mathcal S_{42}$ are pairwise disjoint.
\begin{enumerate}
\item
Assume that $(\mathcal S_3\cup \mathcal S_{42})\cap \mathcal S_2\ne \emptyset$. Then there exists a composition $I=PQR=P'Q'R'$, where
$P$,
$Q$,
$R$,
$P'$,
$Q'$ and $R'$ 
are compositions such that 
$PQR\in \mathcal S_2$,
$P'Q'R'\in \mathcal S_3\cup \mathcal S_{42}$,
\begin{align*}
(\abs{P},\,
\abs{Q},\,
\abs{R})
&=(2m+2,\,
2m+3,\,
2m-1),
\quad\text{and}\\
(\abs{P'},\,
\abs{Q'},\,
\abs{R'})
&=(2m+2,\,
2m+1-r'_{\fo},\,
2m+1+r'_{\fo}).
\end{align*} 
Then $P$ is the prefix of $I$ of size $2m+2$, and so is $P'$. It follows that 
\[
P=P' 
\quad\text{and}\quad
QR=Q'R'.
\] 
Since $P'Q'R'\in\mathcal S_3\cup\mathcal S_{42}$, we have $r'_1=2$. Thus we can write 
\[
I=PQ'2M\!R,
\quad\text{where
$Q'2M=Q$, 
$2M\!R=R'$, and
$\abs{M}=r_{\fo}'$.}
\]
Since $r_{\fo}'$ is odd, $M$ has at least one odd part. Since $R'=2M\!R$, the composition~$M$ contains the first odd part of $R'$, which is $r_{\fo}'$. Thus $M$ consists of the single part $r'_{\fo}$. It follows that the last part of $Q=Q'2M$ is $r_{\fo}'$, which is odd. On the other hand, since $I\in\mathcal S_2$, $q_{-1}$ is even, a contradiction. This proves 
\[
(\mathcal S_3\cup \mathcal S_{42})\cap \mathcal S_2=\emptyset. 
\]
\item
Assume that $\mathcal S_2\cap \mathcal S_{41}\ne\emptyset$. Then there exists a composition 
\[
I=PQR=P'Q'R',
\]
where $PQR\in \mathcal S_2$, $P'Q'R'\in \mathcal S_{41}$,
\begin{align*}
(\abs{P},\,
\abs{Q},\,
\abs{R})
&=(2m+2,\,
2m+3,\,
2m-1),
\quad\text{and}\\
(\abs{P'},\,
\abs{Q'},\,
\abs{R'})
&=(2m+2,\,
2m+1+q'_{\lo},\,
2m+1-q'_{\lo}).
\end{align*}
Then 
\[
P=P' 
\quad\text{and}\quad
QR=Q'R'
\]
as before. Since $q'_{\lo}\ge 3$, we can write 
\[
I=PQM\!R',\quad
\text{where
$QM=Q'$,
$MR'=R$, and
$\abs{M}=q'_{\lo}-2$.}
\]
Since $\abs{M}$ is odd, the composition $M$ has at least one odd part. Since $M$ is a suffix of $Q'$, it contains the last odd part of $Q'$, and $\abs{M}\ge q_{\lo}'$, a contradiction. Thus 
\[
\mathcal S_2\cap \mathcal S_{41}=\emptyset.
\]
\item
Assume that $\mathcal S_3\cap \mathcal S_{42}\ne\emptyset$. Then there exists a composition 
\[
I=PQR=P'Q'R',
\]
where $PQR\in \mathcal S_3$, $P'Q'R'\in\mathcal S_{42}$, and
\begin{align*}
(\abs{P},\,
\abs{Q},\,
\abs{R})
&=(2m+2,\,
2m+1-r_{\fo},\,
2m+1+r_{\fo}),
\quad\text{and}\\
(\abs{P'},\,
\abs{Q'},\,
\abs{R'})
&=(2m+2,\,
2m+1-r_{\fo}',\,
2m+1+r_{\fo}').
\end{align*}
Then
\[
P=P'
\quad\text{and}\quad
QR=Q'R'
\]
as before. Assume that $\abs{Q}\ne\abs{Q'}$. Then we can suppose that $\abs{Q}<\abs{Q'}$ without loss of generality. It follows that $r_{\fo}>r_{\fo}'$. Therefore, the first odd part $r_{\fo}$ of~$R$ is not contained in $R'$, and must lie in~$Q'$. Then 
\[
\abs{Q'}\ge \abs{Q}+r_{\fo}=2m+1,
\]
a contradiction. This proves 
\[
\abs{Q}=\abs{Q'}
\quad\text{and}\quad
r_{\fo}=r_{\fo}'.
\]
Since $QR=Q'R'$, we deduce that $Q=Q'$. It follows that $R=R'$. In summary, we have 
\[
(P,Q,R)=(P',Q',R').
\] 
Since $P'Q'R'\in \mathcal S_{42}$, we have 
\[
\fo(R')-1\le \ell(Q')-\lo(Q')
\quad\text{and}\quad
\fo(R)-1\le \ell(Q)-\lo(Q).
\]
Thus the second alternative in the definition of $\mathcal S_3$ cannot occur.
Hence the first alternative must hold, and so $r_{\fo(R)+1}$ is odd.
This contradicts the condition $r'_{\fo(R')+1}$ is even in $\mathcal S_{42}$.
Therefore, 
\[
\mathcal S_3\cap \mathcal S_{42}=\emptyset. 
\]
\item
Assume that $(\mathcal S_3\cup\mathcal S_{42})\cap \mathcal S_{41}\ne\emptyset$. Then there exists a composition 
\[
I=PQR=P'Q'R',
\]
where $PQR\in \mathcal S_3\cup\mathcal S_{42}$, $P'Q'R'\in\mathcal S_{41}$, 
\begin{align*}
(\abs{P},\,
\abs{Q},\,
\abs{R})
&=(2m+2,\,
2m+1-r_{\fo},\,
2m+1+r_{\fo}), \quad\text{and}\\
(\abs{P'},\,
\abs{Q'},\,
\abs{R'})
&=(2m+2,\,
2m+1+q'_{\lo},\,
2m+1-q'_{\lo}).
\end{align*}
As before, we have 
\[
P=P',\quad 
QR=Q'R',
\]
and we can write
\[
I=PQM\!R',\quad
\text{where 
$QM=Q'$, 
$M\!R'=R$, and
$\abs{M}=r_{\fo}+q'_{\lo}$.}
\]
We proceed according to the number of odd parts of $M$.

Suppose that $M$ has no odd parts. Since $QM=Q'$, we have 
\[
\ell(Q)-\lo(Q)+1\le \ell(Q')-\lo(Q');
\]
since $R=M\!R'$, we have 
\[
\fo(R)-1\ge \fo(R');
\]
since $P'Q'R'\in \mathcal S_{41}$, we have
\[
\fo(R')-1\ge \ell(Q')-\lo(Q').
\]
Therefore,
\begin{equation}\creflabel[ineq]{eq:eQ}
\fo(R)-1\ge \fo(R')\ge \ell(Q')-\lo(Q')+1\ge \ell(Q)-\lo(Q)+2,
\end{equation}
which is impossible since $PQR\in\mathcal S_3\cup\mathcal S_{42}$.

Otherwise, $M$ has at least two odd parts, since $\abs{M}=r_{\fo}+q'_{\lo}$ is even. Thus the first and last odd parts of $M$ are distinct. Since $Q'=QM$, the last odd part of $M$ has value $q_{\lo}'$; since $R=MR'$, the first odd part of $M$ has value $r_{\fo}$. Therefore, the composition $M$ contains two parts: one is of value $r_{\fo}$, the other is of value $q_{\lo}'$. Together with $\abs{M}=r_{\fo}+q_{\lo}'$, we derive that $M=r_{\fo}\,q_{\lo}'$. Since $R=MR'$, we obtain $r_1=r_{\fo}$ is odd. Since $PQR\in\mathcal S_3\cup\mathcal S_{42}$, we have $r_1=2$, a contradiction. This proves 
\[
(\mathcal S_3\cup\mathcal S_{42})\cap \mathcal S_{41}=\emptyset.
\]
\end{enumerate}
In summary, the sets 
$\mathcal S_2$, 
$\mathcal S_3$,
$\mathcal S_{41}$ 
and $\mathcal S_{42}$ are pairwise disjoint. Together with the injectivity of $\varphi_{41}$ and $\varphi_{42}$, the disjointness of $\mathcal S_{41}$ and $\mathcal S_{42}$ shows that the positive terms used for $\mathbf T_4$ are distinct, except for the intentional case $K=L$ handled above by using the common image only once. This completes the proof.
\end{proof}

\section*{Acknowledgments}
The authors are grateful to the anonymous referee for a careful reading of the manuscript and for helpful comments and suggestions. The referee's comments helped us improve the presentation of the motivation in the introduction. This paper was completed while the second author was visiting Professor Jean-Yves Thibon at LIGM, Universit\'e Gustave Eiffel; he thanks Professor Thibon and LIGM for their hospitality.


\end{document}